\documentclass[12pt,oneside, reqno]{amsart}
\usepackage{txfonts}
\usepackage{amsfonts}
\usepackage{mathrsfs}
\usepackage{amsmath}
\usepackage{amssymb, amsmath}
\usepackage{bbm}
\usepackage{stmaryrd}
\usepackage{enumerate}
\newcommand{\be}{\begin{eqnarray}}
\newcommand{\ee}{\end{eqnarray}}
\newcommand{\ce}{\begin{eqnarray*}}
\newcommand{\de}{\end{eqnarray*}}
\newtheorem{theorem}{Theorem}[section]
\newtheorem{lemma}[theorem]{Lemma}
\newtheorem{remark}[theorem]{Remark}
\newtheorem{definition}[theorem]{Definition}
\newtheorem{proposition}[theorem]{Proposition}

\newtheorem{corollary}[theorem]{Corollary}

\def\e{\varepsilon}

\def\a{\alpha}
\def\om{\omega}

\def\p{\partial}
\def\d{\delta}

\def\lam{\lambda}
\def\la{\langle}

\def\[{{\Big[}}
\def\]{{\Big]}}
\def\<{{\langle}}
\def\>{{\rangle}}
\def\({{\Big(}}
\def\){{\Big)}}

\def\bx{{\mathbf{x}}}

\def\min{{\mathord{{\rm min}}}}

\def\no{\nonumber}
\def\bt{\begin{theorem}}
\def\et{\end{theorem}}
\def\bl{\begin{lemma}}
\def\el{\end{lemma}}
\def\br{\begin{remark}}
\def\er{\end{remark}}
\def\bx{\begin{Example}}
\def\ex{\end{Example}}
\def\bd{\begin{definition}}
\def\ed{\end{definition}}
\def\bp{\begin{proposition}}
\def\ep{\end{proposition}}
\def\bc{\begin{corollary}}
\def\ec{\end{corollary}}
\def\cA{{\mathcal A}}

\def\cC{{\mathcal C}}

\def\cF{{\mathcal F}}

\def\cS{{\mathcal S}}

\def\mN{{\mathbb N}}

\def\mR{{\mathbb R}}

\def\sF{{\mathscr F}}

\def\sS{{\mathscr S}}

\def\geq{\geqslant}
\def\leq{\leqslant}

\def\vph{\varphi}
\def\s{\sigma}
\def\bx{{\bf x}}
\allowdisplaybreaks

\allowdisplaybreaks

\begin{document}

\title{Penalization of Reflected SDEs and Neumann Problems of HJB Equations}
\date{}
\author{Jiagang Ren, Jing Wu}

\date{}
\dedicatory{School of Mathematics and Computational Science,
Sun Yat-sen University, \\
Guangzhou, Guangdong 510275, P.R.China\\
Emails: J. Ren: renjg@mail.sysu.edu.cn, J. Wu: wjjosie@hotmail.com}

\begin{abstract}
In this paper we first study the penalization approximation of stochastic differential equations reflected in a domain which satisfies conditions (A) and (B) and prove that the sequence of solutions of the penalizing equations
converges in the uniform topology to the solution of the corresponding reflected stochastic differential equation. Then by using this convergence result, we consider partial differential equations with Neumann boundary conditions
in domains neither smooth nor convex and prove the existence and comparison principle of viscosity solutions of such nonlinear PDEs. Also, by applying the support of reflected diffusions established in \cite{ren-wuAP}, we establish the maximum principle
for the viscosity solutions of linear PDEs with Neumann boundary conditions.
\end{abstract}


\maketitle

\def\px{\pi(x)}
\def\bard{\bar{D}}
\def\lg{\langle}
\def\rg{\rangle}
\def\py{\pi(y)}
\def\nphi{\nabla\vph}
\def\w{\wedge}
\def\v{\vee}
\def\ve{\varepsilon}
\section{Introduction and Background}
 The aim of the paper is to use probabilistic methods to provide the existence, comparison and maximum principle results for viscosity solutions of nonlinear parabolic (and elliptic) equations with Neumann boundary conditions in domains satisfying the exterior sphere condition and weak interior cone property. To this aim, we consider the optimal control problem of diffusion processes governed by the following controlled reflected stochastic differential equation (RSDE for short):
\be\label{eq1}
\begin{cases}
dX(t)=\s(t,X(t), \iota(t))dw(t)+b(t,X(t),\iota(t))dt+d\xi(t), ~~~~t\in[0,T],\\
X(0)=x\in \bard,
\end{cases}
\ee
where $\{w(t)\}_{t\geq 0}$ is a standard $d_1$-dimensional Brownian motion defined on some probability space
$(\Omega,\cF, P)$,  while all together $\nu:=(\Omega,\cF, P,w)$ will be called a reference system, $\sigma$ and $b$ are functions from $[0,T]\times\mR^d\times U$ to $\mR^{d\times d_1}$ and $\mR^d$ respectively, $\iota(\cdot)\in\mathcal{A}_\nu$ and
$\mathcal{A}_\nu$ is the set of progressively measurable processes (w.r.t. the natural filtration generated by $\{w(t)\}$) taking values in a compact metric space $U$. We define the cost function concerning reflected diffusion $X^{t,x}(\cdot)$ starting from $x$ at time $t$:
\be\label{cfV}
V(t,x):=\inf_\nu\inf_{\iota\in\mathcal{A}_\nu}E[\int_t^Tg(s, X^{t,x}(s),\iota(s))ds+h(X^{t,x}(T))].
\ee
Here $\nu$ runs over the set of all the reference systems, $g$ and $h$ are functions defined on $[0,T]\times\mR^d\times U$ and $\mR^d$
respectively. We want to prove, under appropriate conditions, that $V$ is the unique viscosity solution of the following Hamilton-Jacobi-Bellman equation with
Neumann boundary condition:
\be\label{eq2}
\begin{cases}
-\frac{\partial u}{\p t}+H(t,x,u, Du, D^2u)=0,~~~~~~(t,x)\in (0,T)\times D,\\
-\frac{\p u}{\p \mathbf{n}}=0,~~~~~~(t,x)\in (0,T)\times\partial D,\\
u(T,x)=h(x),~~~~~~x\in\bard,
\end{cases}
\ee
where $\mathbf{n}$ are inner normals at $x$ and
$$
H(t,x,r,q,Q):=\max_{\iota\in U}\{-\frac12\mathrm{tr}\sigma\sigma^*(t,x,\iota)Q-\lg b(t,x,\iota),q\rg-g(t,x,\iota)\}.
$$
 Notice that the normals are not necessarily unique since the boundary may have corner points.

In the case of the full space or a domain without reflection, i.e., for diffusions in $\mR^d$ or diffusions stopped upon touching the boundary (hence the Dirichlet boundary), this problem has got systematically studied, and here we refer to P.-L. Lions' pioneering works \cite{lions1, lions2, lions3} and Fleming and Soner's book \cite{fleming}.

If reflection does occur on the boundary, we refer to \cite{ls, sai} for the existence and uniqueness of solutions to reflected stochastic differential equations, where the most general situations (in neither smooth nor convex domains) are considered. In respect of the corresponding elliptic or parabolic equations, a large literature exist, see e.g \cite{bass, frei} (and references therein) for the linear case, and \cite{za13, barles1, barles2, users, ishii1, koi, lions4, za12,ren-wuNA} for the nonlinear case. Yet these results, which are obtained through analytic or probabilistic approaches, are constrained to smooth or at least convex domains. To the best of our knowledge, this problem has never been touched for such general domains as treated in the present paper.

Compared with the case in  smooth domains, in our case there is in the very beginning a notional difficulty as to how to comprehend the Neumann boundary condition at points where the normals are not unique. While in the convex case, the equations can be studied in the context of stochastic variational inequalities by taking advantage of the maximal monotonicity of subdifferentials of the indicator functions of convex sets (see \cite{za12,ren-wuNA}). This decisive tool in the study there is not, unfortunately, appropriate for our general domains.

Hence our first task in this work is to design a penalization scheme to approximate the original RSDE. In the case of reflection in convex domains, Lions, Menaldi and Sznitmann (\cite{lms}) and Slominski in the recent work \cite{sl1} have established the approximation results, where essential estimates and results are obtained through the tool of convex analysis. For domains without smoothness or convexity, however, it
  has been much less well understood. To our best knowledge, so far \cite{saita} is the only one treating the general domains where only the particular case of reflected Brownian motions has been investigated and an almost sure convergence result has been obtained.

  Intuitively, one might image that
 compared with the SDE without reflecting boundary, the main special feature of RSDEs is that a wall is constructed standing vertically along the boundary of $D$ to force the diffusion described by Equation (\ref{eq1})
to remain in $D$. Then the penalization scheme is first to relax the restriction by replacing the vertical wall with
a wall with large slopes (in this way the diffusion is described by an ordinary stochastic differential equation), and then penalize the diffusions by letting the slopes go to infinity. Now a natural question arises: does this penalization sequence of diffusions converge to the original one?
The answer is given in Section 2.

Subsequently, we can ask the following questions: can we understand the equation (\ref{eq2}) via the penalization and is the problem well-posed? We shall give affirmative answers in Sections 3-4.

The third problem we address will be the maximum principle of the viscosity solutions of Neumann problems in``general" domains.
In the seminal paper \cite{str-vara}, Stroock and Varadhan have applied the support of
diffusions generated by usual SDEs to describe the strong maximum principle for corresponding PDEs. For the reflected case, we  presented in the former work \cite{ren-wuAP} the precise characterization of the support of reflected diffusions described by Equation (\ref{eq1}),  and gave  an elementary application
to subharmonic functions (in the viscosity sense) in general domains with reflecting boundary.
 Here in Section 5 we shall apply this result to study more generally the maximum principle for the viscosity solutions of Neumann problems in such domains.

Now let us first specify the assumptions imposed on the domain $D$.

(A) Denoting for $x\in \partial D$ and $r>0$
\ce
&&N_{x,r}:=\{{\bf n}\in \mR^d: |{\bf n}|=1, \mbox{ and } B(x-r{\bf n},r)\cap D=\emptyset\},\\
&&N_x:=\cup_{r>0}N_{x,r},
\de
there exists an $r_0>0$ such that $N_x=N_{x,4r_0}\neq \emptyset$ for all $x\in \partial D$.

(B) There exist constants $\delta>0$ and $\beta\geq 1$ such that for any $x\in \partial D$ there exists a unit vector $l_x$ satisfying
$$
\langle l_x,{\bf n}\rangle \geq\frac{1}{\beta}, ~~\forall {\bf n}\in
\cup_{y\in B(x,\delta)\cap {\partial D}}\cup_{r>0}N_{y,r}.
$$

These two assumptions are presented first in \cite{sai}, where connections with those conditions put forward in \cite{ls} are explained as well.
We notice that an alternate way to express  Condition (A) in nonsmooth analysis is to say that $D$ is
$4r_0$-proximally smooth. That is, putting
$$
D_r:=\{x\in\mR^d, 0<d(x,D)<r\},
$$
where $d(x,D)$ is the distance from $x$ to $D$,
then $d(x,D)$ is $C^{1}$ in $D_{4r_0}$  and the outer boundary of $D_{4r_0}$ is a $C^1$-manifold
(see, e.g. \cite{saita}).

The above two hypotheses (A)-(B) will be in force throughout the paper.

\section{Penalization of reflected SDEs}
\subsection{Preliminary lemmas}
We fix a constant, time horizon, $T>0$ and denote by $\cC([0,T])$ the set of continuous functions defined on $[0,T]$, by $\cC_0([0,T])$ its subset consisting of functions null at zero. For $\a\in (0,\frac12)$ we set
$$
\cC^\a([0,T]):=\{f\in \cC[0,T]:~ |f|_\a<\infty\}
$$
where
$$
|f|_\a:=|f|_0+[f]_\a,
$$
with
$$
|f|_0=\sup_{t\in [0,T]}|f(t)|;~~[f]_\a:=\sup_{s,t\in [0,T], s\neq t}
\frac{|f(s)-f(t)|}{|s-t|^\a}.
$$
We also set $\cC_0^\a([0,T]):=\cC_0([0,T])\cap \cC^\a([0,T])$. If $0\leq s\leq t\leq T$, we define
$$
[f]_{\a;s,t}:=\sup_{u,v\in [s,t], u\neq v}
\frac{|f(u)-f(v)|}{|u-v|^\a};
$$
$$
|f|_{0;s,t}:=\sup_{s\leq u\leq t}|f(u)|;
$$
$$
\Delta_{s,t}(f):=\sup_{s\leq u\leq v\leq t}|f(v)-f(u)|,
$$
$$
|f|^s_t:=\sup\sum_{k=1}^n|f(t_k)-f(t_{k-1})|,
$$
where the supermum is taken over all partitions $s=t_1<t_2<\cdots t_n=t$ and $n=1,2,3,\cdots$.

Let $\rho:\mR_+\to\mR_+$ be a $\cC^1$ increasing function such that
$$
\rho(t)=
\begin{cases}
t, & t\in [0,4r^2_0];\\
9r_0^2, & t\geq 9r^2_0.
\end{cases}
$$

Set $\vph_0(x)=d(x,\bard)^2$ and $\vph(x)=\rho(\vph_0(x))$.
The function $\vph$ has the following properties.

\bl\label{lemma1}
(i) $\vph\in \cC^{1,1}_b(\mR^d)$, where $\cC^{1,1}_b(\mR^d)$ stands for the space of bounded functions on $\mR^d$ with bounded and Lipschitz continuous first derivatives.

(ii) Every $x\in D_{4r_0}$ has a unique projection $\pi(x)$ on $\bard$, and
$\vph(x)=|x-\px|^2$ for $x\in D_{2r_0}$.

(iii) $\vph(x)=$constant for $x\notin D_{3r_0}$.

(iv) For $n\in \mN_+$, $x,y\in D_{2r_0}$ we have
\ce
&&\lg \nabla\vph(x), x-y\rg\geq
-2\rho'(\vph_0(x))\{\lg y-\py,x-\px\rg+\frac{1}{8r_0}|\px-\py|^2|x-\px|\}.
\de

(v) Let $x\in D_{2r_0}$. If either $y\in\partial D$, $z\in N_y$
or $y\in D$, $z=0$, then \ce
&&\lg x-y, z\rg\geq
\lg x-\px, z\rg-\frac{1}{8r_0}|\px-y|^2.
\de

\el

\begin{proof}

(i)-(iii) are proven
\cite[Section 1]{saita}.

It is easy to see that Condition (A) implies that
for $x\in \partial D$ and $z\in N_{x}$,
\be\label{Ain}
\lg y-x, z\rg \geq -\frac{1}{8r_0}|y-x|^2, ~\forall y\in \bard.
\ee

Thus for $y\in\partial D$ and $z\in N_y$, $x\in D_{2r_0}$,
\ce
\lg x-y,z\rg&=&\lg x-\pi(x),z\rg+\lg \pi(x)-y,z\rg\\
&\geq&\lg x-\pi(x),z\rg-\frac{1}{8r_0}|\px-y|^2,
\de
and (v) is proved.

To prove (iv), note that for $x,y\in D_{4r_0}$ and
$$
z\begin{cases}
\in N_{\pi(x),4r_0},~&\mbox{if}~x\notin D;\\
=0,~&\mbox{if}~x\in D,
\end{cases}
$$
we have by (\ref{Ain}),
\ce
\lg y-x,z\rg&=&\lg y-\pi(y),z\rg+\lg\pi(y)-\pi(x),z\rg+\lg\pi(x)-x,z\rg\\
&\geq&\lg y-\pi(y),z\rg-\frac{1}{8r_0}|\pi(y)-\pi(x)|^2.
\de

Thus for $x,y\in D_{2r_0}$,
\ce\label{inequality2}
&&\lg x-y, \nabla \vph(x)\rg\no\\
&=&2\lg x-y, \rho'(\vph_0(x))(x-\pi(x))\rg\no\\
&\geq&
-2\rho'(\vph_0(x))\big[\lg y-\py, x-\pi(x)\rg+\frac{1}{8r_0}
|\py-\px|^2|x-\pi(x)|\big].
\de
\end{proof}

For $z\in \cC^\a([0,\infty))$, an $\a$-H\"older continuous function on $[0,\infty)$, $x\in \bard$ and positive integer $n\in \mN_+$, consider the following deterministic differential equation
$$
y_n(t)=x+z(t)-\frac{n}{2}\int_0^t\nphi(y_n(s))ds.
$$

Let
\ce\label{ve}
\ve(s,t):=
\frac{\delta}{8(1+4\beta+\beta^2\exp\{\frac{\beta l^2}{r_0}(|z|_{0;s,t}+\delta)\})}
\wedge r_0,
\de
\ce
\ve(t):=\ve(0,t),
\de
where $l$ is the Lipschitz constant of the projection $\pi(\cdot)$,
and
$$
\xi(t):=\frac{n}{2}\int_0^t\nphi(y_n(s))ds,
$$
and we use $\xi^\tau$ to denote the stopped function $t\mapsto \xi(t\wedge \tau)$
for a stopping time $\tau$.

We will need the following result.
\bl\label{lem22}
Set
$$
\tau^n:=\inf\{v: \sup_{0\leq u\leq v}d(y_n(u),\overline{D})\geq \ve(v)\}\wedge T.
$$
Then there exist a constant $C\geq0$
such that for $0\leq s<t\leq T$ we have for $\a\in(0,1)$,
$$
|\xi^{\tau^n}|^s_t\leq C\Big(1+(t-s)\d^\frac{1}{\a}|z^{\tau^n}|^\frac{1}{\a}_{\a;s,t}\exp\{C(1+|z^{\tau^n}|_{\a;s,t})
\}\Big)\{\Delta_{s,t}(z^{\tau^n})+1\}.
$$
\el

\begin{proof}
We use an argument similar to the proof of \cite[Prop. 5.1]{saita}.
Since $n$ will be fixed in the proof, we will omit it for notational
simplicity.

Define inductively
$$
\tau_0:=s;
$$
$$
\tau_{2k+1}=\inf\{u> \tau_{2k}: |\pi(y_n(u))-\pi(y_n(\tau_{2k}))|\geq \frac{\delta}{2}\}\wedge t,~k\geq 0,
$$
$$
\tau_{2k}=\inf\{u\geq \tau_{2k-1}: \pi(y_n(u))\in \partial D\}\wedge t,~k\geq 1.
$$
Then by \cite[Lemma 5.1]{saita} we have for $\tau_{2k}\leq u<v\leq \tau_{2k+2}$
\be
|\xi^\tau|^u_v\leq \beta(\Delta_{u,v}(y_n^\tau)+\Delta_{u,v}(z^\tau));
\ee
and by \cite[Lemma 5.3]{saita},
$$
\Delta_{u,v}(y_n^\tau)\leq K(u,v)(\Delta_{u,v}(z^\tau)+\ve(u,v)),
$$
where
$$K(u,v)=8\beta\exp\{\frac{\beta l^2}{2r_0}(|z^\tau|_{0;s,t}+\delta)\}+1,$$
where $l$ is the Lipschitz constant of $\pi$.

Hence
\be
|\xi^\tau|^u_v\leq K_2(u,v)\{\Delta_{u,v}(z^\tau)+\ve(s,t)\},
\ee
where $K_2(u,v)=2\beta+8\beta^2\exp\{\frac{\beta l^2}{2r_0}(|z^\tau|_{0;s,t}+\delta)\}$.
Since
\ce
&&|y_n^\tau(\tau_{2k+1})-y_n^\tau(\tau_{2k})|\\
&\geq& |\pi(y_n^\tau(\tau_{2k+1}))-\pi(y_n^\tau(\tau_{2k}))|-2\ve(u,v)\\
&\geq&\frac{\delta}{2}-2\ve(u,v),
\de
we have
$$
\frac{\delta}{2}-2\ve(u,v)\leq K(u,v)
\{\Delta_{\tau_{2k}, \tau_{2k+1}}(z^\tau)+\ve(u,v)\}.
$$
Hence
\ce
\Delta_{\tau_{2k}, \tau_{2k+1}}(z^\tau)&\geq&
\frac{1}{K(\tau_{2k},\tau_{2k+1})}\big(\frac{\delta}{2}-2\ve(\tau_{2k},\tau_{2k+1})\big)-
\ve(\tau_{2k},\tau_{2k+1})\\
&\geq&\frac{\delta-4\ve(\tau_{2k},\tau_{2k+1})}{2K(\tau_{2k},\tau_{2k+1})}
-\ve(\tau_{2k},\tau_{2k+1})\\
&\geq& \frac{\delta}{4K(\tau_{2k},\tau_{2k+1})},
\de
where in the last step we used the fact that
$$
\ve(\tau_{2k},\tau_{2k+1})\leq \frac{\delta}{4(K(\tau_{2k},\tau_{2k+1})
+2)}.
$$
Thus
\ce
\tau_{2(k+1)}-\tau_{2k}\geq \left[\frac{\delta}{4K(\tau_{2k},\tau_{2k+1})|z^\tau|_{\a; \tau_{2k}, \tau_{2(k+1)}}}\right]^\frac{1}{\a}.
\de
Letting $N$ be the smallest $k$ such that $\tau_{2k}=t$ we have
\ce
t-s\geq \sum_{k=0}^{N-1}(\tau_{2(k+1)}-\tau_{2k})\geq (N-1)\left[\frac{\delta}{4K(s,t)|z^\tau|_{\a; s,t}}\right]^\frac{1}{\a}.
\de
Therefore
\ce
N\leq (t-s)\left[\frac{4K(s,t)|z^\tau|_{\a; s, t}}{\delta}\right]^\frac{1}{\a}+1.
\de
Substituting this into the last expression of the following
inequalities we complete the proof:
\be
|\xi^\tau|^s_t &\leq& \sum_{k=0}^N|\xi^\tau|^{\tau_{2k}}_{\tau_{2(k+1)}}\nonumber\\
&\leq&NK_2(s,t)\{\Delta_{s,t}(z^\tau)+\ve(s,t)\}\nonumber\\
&\leq&CNK_2(s,t)\{\Delta_{s,t}(z^\tau)+1\}.
\ee
\end{proof}

The following simple result, which may be considered as a stochastic Gronwall lemma, will play
an important role in the sequel.

\bl\label{inequality1}
Let $\{M_t\}$ be a local martingale, $\{G_t\}$ and $\{F_t\}$ be increasing processes with $F_0=G_0=0$, and $\{H_t\}$ be a nonnegative
semimartingale. If
$$
dH_t\leq dM_t+dG_t+H_tdF_t,~~\forall t,
$$
then for any bounded stopping time $\tau$ we have
\ce
E[H_\tau e^{-F_\tau}]\leq E[\int_0^\tau e^{-F_t}dG_t]
\de
and
\ce
E[\sup_{t\leq \tau}H_t^\alpha e^{-\alpha F_t}]\leq \frac{2-\alpha}{1-\alpha}E[\int_0^\tau e^{-F_t}dG_t]^\alpha, ~~\alpha\in (0,1).
\de
\el
\begin{proof}It suffices to prove the first inequality since then the second one will follow by \cite[Ch. IV, Prop. 4.7]{ry}.
Let
$$
N_t:=H_te^{-F_t}.
$$
By Ito formula we have
\ce
dN_t\leq e^{-F_t}dM_t+e^{-F_t}dG_t.
\de
Hence, as desired,
$$
E(N_\tau)\leq E\big[\int_0^\tau e^{-F_t}dG_t\big].
$$
\end{proof}

\subsection{Penalization results}
We denote by $\{\cF_t\}$ the natural filtration generated by $\{w(t)\}$.
Then a solution of Equation (\ref{eq1}) is a pair of $\{\cF_t\}$-adapted continuous processes
$(X(t),\xi(t))$ such that

(i) $(X(t))$ is $\bard$-valued;

(ii) $(\xi(t))$ is of locally bounded variation and $\xi_0=0$;

(iii) $(X(t),\xi(t))$ satisfies the equation (\ref{eq1}) and
$$
\xi(t)=\int_0^t{\bf n}(s)d|\xi|(s),
$$
$$
|\xi|(t)=\int_0^t \mathbbm{l}_{\partial D}(X(s))d|\xi|(s),
$$
where ${\bf n}(s)\in N_{X(s)}$ if $X(s)\in\partial D$ and $|\xi|(\cdot)$
is the total variation process of $\xi(\cdot)$.

(H) $\sigma$ and $b$ are bounded functions, Lipschitz continuous with respect to $x$, and the supports of $\sigma$ and $b$ are contained in
$[0,T]\times(D_{3r_0}\cup\bard)\times U$.

Note that for the reflected equation (\ref{eq1}) itself, what do matter are the values of $\sigma$ and $b$ on $\bard$.
Nevertheless, since every $x\in D_{4r_0}$ has a projection on $\bard$, one can easily extend a function $g$ defined
only on $\bard$ to $D_{4r_0}\cup\bard$ by letting $g(x):=g(\pi(x))$ and then one can modify it in a stand way to make it
supported by $D_{3r_0}\cup\bard$.

It is proved in \cite{sai} that under the assumptions (A)-(B) and (H) Equation (\ref{eq1}) has a unique strong solution.

The penalization scheme considered is defined by the sequence of SDEs
\be\label{pn}
\begin{cases}
dX_n(t)=\s(t,X_n(t),\iota(t))dw(t)+b(t,X_n(t),\iota(t))dt-\frac{n}{2}\nphi(X_n(t))dt\\
X_n(0)=x\in \bard
\end{cases}
\ee

Since $\nphi$ has compact support and is Lipschitz, it is classic that this equation has a unique solution (denoted by $(X_n,\xi_n)$ for simplicity).

For a continuous function $f: [0,T]\mapsto \mR$ we define for $h>0$
$$
\om_h(f):=\max_{s,t\in [0,T], |s-t|<h}|f(s)-f(t)|.
$$

The following result is proved in \cite{sl1} when $D$ is convex and $U$ is a singleton.

\bl\label{lem32}
$$
\sup_{\iota}\|\sup_{0\leq t\leq T}d(X_n(t),\bard)\|_p\leq C(\frac{\ln n}{n})^\frac12.
$$
$$
\sup_{\iota}\|\sup_{k\geq n}\sup_{0\leq t\leq T}d(X_k(t),\bard)\|_p\leq
C\frac{(\ln n)^\frac{1}{2}}{n^{\frac12-\frac{1}{p}}}, ~~\forall p>2.
$$
\el

\begin{proof}
Assume for simplicity $T=1$ and let $2K$ be the Lipschitz constant of $\nphi$ and set
$$
t_n^k:=\frac{k}{n},~~k=0,1,\cdots, n.
$$
Set
$$
Y^\iota_n(t):=\int_0^t\s(s,X_n(s),\iota(s))dw(s)+\int_0^tb(s,X_n(s),\iota(s))ds.
$$
Then for $t\in [t_n^k, t_n^{k+1}]$,
$$
X_n(t)=X_n(t_n^k)+Y^\iota_n(t)-Y^\iota_n(t_n^k)-\frac{n}{2}\int_{t_n^k}^t\nphi(X_n(s))ds.
$$
Now we fix $\om\in A_n$ where
$$
A_n:=\{\om_{\frac{1}{n}}(Y^\iota_n)e^K\leq 2r_0\},
$$
and omit the variable $n$ (as always).

For $t\in [t_n^0, t_n^1]$, set
$$
Z_n(t)=x.
$$
Then $Z_n$ satisfies the equation
$$
Z_n(t)=x-\frac{n}{2}\int_0^t\nphi(Z_n(s))ds,~t\in [t_n^0, t_n^1].
$$
Hence
$$
|X_n(t)-Z_n(t)|\leq |Y^\iota_n(t)|+Kn\int_0^t|X_n(s)-Z_n(s)|ds.
$$
We deduce by Gronwall lemma
$$
|X_n(t)-Z_n(t)|\leq \om_{\frac{1}{n}}(Y^\iota_n)e^K,~~t\in [t_n^0, t_n^1].
$$
In particular, $d(X_n(t_n^1),\bard)\leq 4r_0$.
Hence $\pi(X_n(t_n^1))$ makes sense.

For $t\in [t_n^1, t_n^2]$, set
$$
Z_n(t)=\pi(X_n(t_n^1))+(X_n(t_n^1)-\pi(X_n(t_n^1)))e^{-n(t-t_n^1)}.
$$
Then $Z_n(t)$ satisfies
$$
Z_n(t)=X_n(t_n^1)-\frac{n}{2}\int_{t_n^1}^t\nphi(Z_n(s))ds.
$$
Then
$$
d(Z_n(t_n^2), \bard)=d(X_n(t_n^1),\bard)e^{-1}\leq \om_{\frac{1}{n}}(Y^\iota_n)e^Ke^{-1},
$$
and, again by Gronwall inequality as we just did,
$$
|Z_n(t)-X_n(t)|\leq \om_{\frac{1}{n}}(Y^\iota_n)e^K
$$
Thus
$$
d(X_n(t_n^2),\bard)\leq \om_{\frac{1}{n}}(Y^\iota_n)e^K
+d(X_n(t_n^1),\bard)e^{-1}.
$$
Continuing this process we obtain that for any $k$
$$
d(X_n(t_n^k),\bard)\leq \om_{\frac{1}{n}}(Y^\iota_n)e^K\sum_{i=0}^\infty e^{-i}\leq 2\om_{\frac{1}{n}}(Y^\iota_n)e^K.
$$

For $t\in [t_n^k, t_n^{k+1}]$,
\ce
Z_n(t)&=&X_n(t_n^k)-\frac n2\int_{t_n^k}^t\nphi(Z_n(s))ds,\\
X_n(t)&=&X_n(t_n^k)+Y^\iota_n(t)-Y^\iota_n(t_n^k)
-\frac{n}{2}\int_{t_n^k}^t\nabla\vph(X_n(s))ds.
\de
By simple calculations we have
\ce
|X_n(t)-Z_n(t)|\leq \om_{\frac{1}{n}}(Y^\iota_n)e^K.
\de

Moreover, since \ce
Z_n(t)=\pi(X_n(t_n^k))+(X_n(t_n^k)-\pi(X_n(t_n^k)))e^{-n(t-t_n^k)},
\de
\ce
d(Z_n(t),\bard)\leq d(X_n(t_n^k),\bard)e^{-1}.
\de
Thus for $t\in [t_n^k, t_n^{k+1}]$,
\ce
\big(d(X_n(t),\bard)&\leq&|X_n(t)-Z_n(t)|+\big(d(Z_n(t),\bard)\\
&\leq&\om_{\frac{1}{n}}(Y^\iota_n)e^K+d(X_n(t_n^k),\bard)e^{-1}\\
&\leq&\om_{\frac{1}{n}}(Y^\iota_n)e^K(1+2e^{-1}).
\de

Summing up we have
$$
\sup_{t\in [0,T]}d(X_n(t),\bard)\leq (1+2e^{-1})\om_{\frac{1}{n}}(Y^\iota_n)e^K.
$$
Remember that the above reasoning is done for $\om\in A_n$. For
$\om\notin A_n$ we use the fact that
$$
\sup_{t\in [0,T]}d(X_n(t),\bard)\leq 1+3r_0,
$$
which follows from the fact that $\sigma$, $b$ and $\nabla \varphi$ are all supported by $D_{3r_0}\cup\bard$.

Then we have
$$
\sup_{t\in [0,T]}d(X_n(t),\bard)\leq C1_{A_n}\om_{\frac{1}{n}}(Y^\iota_n)
+1_{A_n^c}C.
$$
It follows from Chebyshev inequality that
$$
\|1_{A_n^c}\|_p\leq C\|\om_{\frac{1}{n}}(Y^\iota_n)\|_p.
$$
Consequently,
$$
\|\sup_{t\in [0,T]}d(X_n(t),\bard)\|_p\leq C\|\om_{\frac{1}{n}}(Y^\iota_n)\|_p,
$$
where according to
\cite[Lemma A.4]{sl}, the RHS is dominated by $C(\frac{\ln n}{n})^\frac12$, with $C$ independent of $\iota$. This proves the first inequality.

For the second one, one only needs to note that for every $\iota$,
\ce
E[|\sup_{k\geq n}\sup_{t\in [0,T]}d(X_k(t),\bard)|^p]
&\leq & \sum_{k=n}^\infty E[|\sup_{t\in [0,T]}d(X_k(t),\bard)|^p]\\
&\leq & C\sum_{k=n}^\infty\frac{(\ln k)^{p/2}}{k^{p/2}}\\
&\leq &C\int_n^\infty x^{-p/2}(\ln x)^{p/2}dx\\
&\leq & Cn^{1-p/2}(\ln n)^{p/2},
\de
where the constants are independent of $\iota$ and the proof is thus finished.
\end{proof}
\bt\label{thm32}
\ce
\lim_{m, n\to\infty}\sup_{\iota\in\cA_\nu} E[\sup_{0\leq s\leq T}|X_n(s)-X_m(s)|^p]=0.
\de
\et
\begin{proof}
Set
$$
\tau_n:=\inf\{t: \sup_{0\leq s\leq t}d(X_n(s),\bard)\geq \ve(t)\}.
$$
$$
Y_{n,m}(t)=X_n(t)-X_m(t),
$$
$$
Y_{n,m}^{\tau_n\w\tau_m}(t)=Y_{n,m}(t\w\tau_n\w\tau_m),
$$
where $\ve(t)$ is defined in (\ref{ve}).

Then we have by It\^o formula
\be\label{for8}
|Y_{n,m}^{\tau_n\w\tau_m}(t)|^2&=&2\int_0^{t\w\tau_n\w\tau_m}(X_n(s)-X_m(s))
(\s(s,X_n(s),\iota(s))-\s(s,X_m(s)),\iota(s))dw(s)\no\\
&&+2\int_0^{t\w\tau_n\w\tau_m}(X_n(s)-X_m(s))
(b(s,X_n(s),\iota(s))-b(s,X_m(s),\iota(s)))ds\no\\
&&-\int_0^{t\w\tau_n\w\tau_m}
(n\nphi(X_n(s))-m\nphi(X_m(s)))(X_n(s)-X_m(s))ds\no\\
&\leq & M_{n,m}(t)+\int_0^{t\w\tau_n\w\tau_m}g_{n,m}(s)ds+\int_0^{t\w\tau_n\w\tau_m}f_{n,m}(s)|Y_{n,m}^{\tau_n\w\tau_m}(s)|^2ds,
\ee
where
\ce
&&f_{n,m}(s)=C(1+n|\nphi(X_n(s))|+m|\nphi(X_m(s))|),\\
&&g_{n,m}(s)=C(n|\nphi(X_n(s))d(X_m(s),\bard)|
+m|\nphi(X_m(s))d(X_n(s),\bard)|),
\de
and $M_\cdot$ is a martingale, and (\ref{for8}) holds since by (iv) in Lemma \ref{lemma1},
\ce
&&-\int_0^{t\w\tau_n\w\tau_m}
(n\nphi(X_n(s))-m\nphi(X_m(s)))(X_n(s)-X_m(s))ds\\
&=&-\int_0^{t\w\tau_n\w\tau_m}
n\nphi(X_n(s))(X_n(s)-X_m(s))ds-\int_0^{t\w\tau_n\w\tau_m}
n\nphi(X_m(s))(X_m(s)-X_n(s))ds\\
&\leq&2\int_0^{t\w\tau_n\w\tau_m}
n\nphi(X_n(s))\Big[\lg d(X_m(s),\bard),d(X_n(s),\bard)\rg-\frac{1}{8r_0}|\pi(X_n(s))-\pi(X_m(s))|^2d(X_n(s),\bard)\Big]\\
&&+2\int_0^{t\w\tau_n\w\tau_m}
m\nphi(X_m(s))\Big[\lg d(X_n(s),\bard),d(X_m(s),\bard)\rg-\frac{1}{8r_0}|\pi(X_m(s))-\pi(X_n(s))|^2d(X_m(s),\bard)\Big]\\
&\leq&\int_0^{t\w\tau_n\w\tau_m}C(n|\nphi(X_n(s))|+m|\nphi(X_m(s))|)|Y_{n,m}^{\tau_n\w\tau_m}(s)|^2ds
+\int_0^{t\w\tau_n\w\tau_m}g_{n,m}(s)ds.
\de
Therefore by Lemma \ref{lem22}-\ref{inequality1} we have
\ce
&&E[\sup_{t\leq T}|Y^{\tau_n\w\tau_m}_{n,m}(t)\exp\{-\frac12\int_0^{t\w\tau_n\w\tau_m}f_{n,m}(s)ds\}]\\
&\leq &3E[\int_0^{T\w\tau_n\w\tau_m}\exp\{-\frac12\int_0^tf_{n,m}(s)ds\}g_{n,m}(t)dt]^\frac12\\
&\leq & 3\left\{E\int_0^{T\w\tau_n\w\tau_m}g_{n,m}(t)dt\right\}^\frac12\\
&\leq &CE[\sup_{t\leq T}d(X_m(t),\bard)^\frac12(\int_0^{T\w\tau_n\w\tau_m}|n\nphi(X_n(s))|ds)^\frac12]\\
&&+E[\sup_{t\leq T}d(X_n(t),\bard)^\frac12(\int_0^{T\w\tau_n\w\tau_m}|m\nphi(X_m(s))|ds)^\frac12]\\
&\leq& C\{E[\sup_{t\leq T}(d(X_m(t),\bard)+d(X_n(t),\bard)]\}^\frac12\\
&&\cdot\{E[\int_0^{T\w\tau_n\w\tau_m}(|n\nphi(X_n(s))|+|m\nphi(X_m(s))|ds]\}^\frac12\\
&\leq& C(n^{-\frac14}(\ln n)^\frac14+m^{-\frac14}(\ln m)^\frac14),
\de
where the last inequality follows from Lemma \ref{lem22} and Lemma  \ref{lem32}.

Hence uniformly in $\iota$,
$$\sup_{t\leq T}|Y^{\tau_n\w\tau_m}_{n,m}(t)\exp\{-\frac12\int_0^{t\w\tau_n\w\tau_m}f_{n,m}(s)ds\}|
\to 0~~\mbox{in probability}.
$$
Next set
\ce
Z_n^\iota=\int_0^t\s(s,X_n(s),\iota(s))dw(s)+\int_0^tb(s,X_n(s),\iota(s))ds.
\de
By  Lemma \ref{lem22} we have that there exists two positive constants $c_1$ and $c_2$ such that
\ce
\int_0^{T\wedge\tau_n\wedge\tau_m}f_{n,m}(s)ds\leq c_1\exp(c_2|Z_n^\iota|_{\a;0,T}).
\de
Let $B$ be the DDS-Brownian motion of the martingale part of $Z_n^\iota$, we have by using the boundedness of $\s$
\ce
|Z^\iota_n|_{\a;0,t}\leq C[|B|_{\a;0,c_3t}+1],
\de
which combined with the Fernique theorem gives
$$
\sup_{n,m}\sup_{\iota\in\cA_\nu} E[\int_0^{T\w\tau_n\w\tau_m}f_{n,m}(s)ds]<\infty.
$$
Consequently $\{\int_0^{T\w\tau_n\w\tau_m}f_{n,m}(s)ds\}_{n,m}$ is bounded in probability uniformly in $\iota, n$ and $m$,
and this in turn implies that
$$
\{\exp\{-\frac12\int_0^Tf_{n,m}(s)ds\}\}_{n,m}
$$
is uniformly bounded away from zero in probability. This yields that
$$\sup_{t\leq T}|Y^{\tau_n\w\tau_m}_{n,m}(t)|
\to 0~~\mbox{in probability uniformly in }\a.
$$
Now, since by Lemma \ref{lem32}
$$
\sup_{\iota\in\cA_\nu} P(\tau_n\w\tau_m<T)\to 0, \quad\mbox{as}\quad n, m\to\infty,
$$
we have in fact
$$\sup_{t\leq T}|Y_{n,m}(t)|
\to 0~~\mbox{in probability}
$$
uniformly in $\iota\in\cA_\nu$. This combined with the fact that
$$
\sup_n\sup_{\iota\in\cA_\nu} E[\sup_{t\leq T}|X_n(t)|^p]<\infty,~~\forall ~p\geq1,
$$
of course implies that
$$
\sup_{\iota\in\cA_\nu} E[\sup_{t\leq T}|X_n(t)-X_m(t)|^p]\to 0, ~~\forall ~p\geq1.
$$
\end{proof}

Now we can prove

\bt\label{pconvergence} For any $p\geq 1$,
\be
\lim_{n\to\infty}\sup_{\iota\in\cA_\nu} E[\sup_{0\leq s\leq T}|X_n(s)-X(s)|^p]=0.
\ee
\et
\begin{proof}
Use the same $\tau_n$ as above and set
$$
Y_n(t)=X_n(t\w\tau_n)-X(t\w\tau_n).
$$
By a similar calculus as above we have
\ce
|Y_n(t)|^2\leq M_n(t)+C\int_0^{t\w\tau_n}f_n(s)|Y_n(s)|^2ds
+C\int_0^{t\w\tau_n}|Y_n(s)|^2d|\xi|_s+G_n(t),
\de
where
\ce
&&M_n(t):=2\int_0^{t\w\tau_n}\lg X_n(s)-X(s),\sigma(X_n(s),\iota(s))-\sigma(X(s),\iota(s))\rg dw(s),\\
&&f_n(t):=C(1+n|\nphi(X_n(t))|),\\
&&G_n(t):=C\int_0^{t\w\tau_n}|d(X_n(s),\bard)d|\xi|_s\leq C\ve(t)|\xi|_t.
\de

The conclusion now follows from similar arguments as in the proof of Theorem \ref{thm32}.
\end{proof}

\section{Cauchy-Neumann boundary problems in $D$}

Suppose that $D$ is a domain satisfying the Conditions (A) and (B) in Section 2.
In this section we consider the following PDE with Neumann boundary condition:
\be\label{pde}
\left\{
\begin{array}{lllll}
-\frac{\partial u}{\partial t}+H(t,x,u,Du,D^2u)=0\quad \mbox{in}~~~(0,T)\times D,\\
\\
-\frac{\partial u}{\partial \mathbf{n}}=0\quad \mbox{on}~~~(0,T)\times \partial D,\\
\\
u(T,\cdot)=h(\cdot)\quad \mbox{on}~~~\bard,
\end{array}\right.
\ee
where $H(t,x,r,q,Q)$ is a continuous function from $[0,T]\times\bard\times\mR\times \mR^d\times \cS^d$ to $\mR$,
$\cS^d$ is the set of all real-valued $d\times d$ symmetric matrices. We assume that $F$ satisfies the following
hypotheses:

(H1)
If $u\leq v$, then there exists a constant $\gamma>0$ such that
\be
\gamma(v-u)\leq H(t,x,v,q,X)- H(t,x,u,q,X),~~~\forall (t,x,q,X)\in[0,T]\times \bard\times\mR^d\times \cS^d;
\ee

(H2) there exists an increasing continuous function $\omega_1: [0,+\infty)\mapsto [0,+\infty)$, $\omega(0)=0$, such that
\be
|H(t,x,u,q_1,X)-H(t,x,u,q_2,Y)|\leq \om_1(|q_1-q_2|+\|X-Y\|)
\ee
for all $(t,x,u)\in (0,T)\times \bard\times\mR$, $q_1,q_2\in \mR^d$, $X,Y\in \cS^d$.

(H3) there exists an increasing continuous function $\omega_2: [0,+\infty)\mapsto [0,+\infty)$, $\omega(0)=0$, such that
\be
H(t,y,u,\a(x-y),Y)-H(t,x,u,\a(x-y),X)\leq \omega_2(|x-y|(1+\alpha|x-y|))
\ee
for all $x,y\in D$, $u\in\mR$, $p\in\mR^d$, $X,Y\in \cS^d$ and  $\a>1$ satisfying
\ce
\label{matrix}
-3\alpha\left(
\begin{array}{llll}
I\quad \quad0
\\
0\quad \quad I
\end{array}\right)\leq  \left(
\begin{array}{llll}
{X}\quad\quad 0
\\
0\quad -{Y}
\end{array}\right)\leq
3\alpha\left(
\begin{array}{llll}
I\quad \quad-I
\\
-I\quad \quad I
\end{array}\right).
\de

Note that (H3) implies the ellipticity of $H$ (c.f. \cite[Proposition 3.8]{koi}).
We denote by $\cC^{1,2}((0,T)\times D)$ the set of all functions $u$ defined in  $(0,T)\times D$ such that
$u$ and all the derivatives $\frac{\partial u}{\partial t}$ and $D^k u$, $|k|\leq 2$, are continuous and bounded in
$(0,T)\times D$. Here $k=(k_1,\cdots,k_d)$ is a multi-index with $k_i=0,1$ and
$$
|k|=k_1+\cdots+k_d,
$$
$$
D^k=D^{k_1}_{x_1}\cdots D^{k_d}_{x_d}.
$$
$C^{1,2}([0,T]\times \bar{D})$ will be the subset of $C^{1,2}((0,T)\times D)$ consisting of such $u$ that all the above involved
derivatives extends to bounded continuous functions on $[0,T]\times\bar{D}$.

In addition to (A) and (B), we also assume

(C) There exists a function $f\in\mathcal{C}_b^2(\mR^d)$ with $\sup_{x\in\mR^d}|f(x)|\leq 1$ and a constant $~\gamma>0$ such that $\forall x\in\partial D$ and $\forall \mathbf{n}\in N_x$,
$$
D f(x)\cdot \mathbf{n}\geq \frac{\gamma}{8r_0},
$$
where $r_0$ is the same as in condition (A).

We define, in a standard way following e.g. \cite{koi, users}, for a function $u:(0,T)\times \bard\mapsto \mR$, the superjet
and subjet

 \ce
 &&J^{1,2,+}u(t,x):=\Big\{(p,q,Q)\in\mR\times\mR^d\times\mathcal{S}^d;\\
~~~~~~~~~~~&&\limsup_{(0,T)\times\bard\ni(s,y)\to(t,x)}\frac{u(s,y)-u(t,x)-p(s-t)-\lg q,y-x\rg-\frac12\la Q(y-x),y-x\rg}{|s-t|+|x-y|^2}\leq0\Big\},\\
 &&J^{1,2,-}u(t,x):=\Big\{(p,q,Q)\in\mR\times\mR^d\times\mathcal{S}^d;\\
 ~~~~~~~~~~&&\liminf_{(0,T)\times\bard\ni(s,y)\to(t,x)}\frac{u(s,y)-u(t,x)-p(s-t)-\lg q,y-x\rg-\frac12\la Q(y-x),y-x\rg}{|s-t|+|x-y|^2}\geq0\Big\}.
 \de
and their closures
\ce
&&\bar{J}^{1,2,\pm}u(t,x):=\Big\{(p,q,Q)\in\mR\times\mR^d\times\mathcal{S}^d;
\exists (t_k,x_k)\in (0,T)\times \bard \mbox{ and } (p_k,q_k,X_k)\in  J^{1,2,\pm}u(t,x) \\
~~~~~~~~~~&&\mbox{ such that } (t_k,x_k, u(t_k,x_k),p_k, q_k, X_k)\to (t, x, u(t,x), p,q, X)\Big\}.
\de

Define for $x\in\partial D$,
\ce
&&N^-_x(q):=\inf_{\mathbf{n}\in N_x}\lg q, -\mathbf{n} \rg,\\
&&N^+_x(q):=\sup_{\mathbf{n}\in N_x}\lg q, -\mathbf{n} \rg.
\de
Since $\vph\in C_b^{1,1}$ and
$$
N_x=\{\frac12\nabla\vph(y)|y-x|^{-1}: y\in D_{2r_0}, \pi(y)=x\},
$$
we have
$$
N_x=\cap_{\delta>0}\cup_{|x'-x|<\delta,x'\in\partial D}\{\frac12\nabla\vph(y)|y-x|^{-1}: y\in D_{2r_0}\cup\partial D, \pi(y)=x'\},
$$
from which we deduce that
\be\label{normals1}
N^-_x(q)=\lim_{\delta\downarrow 0}\inf\{\lg q, -\mathbf{n'} \rg: |x'-x|<\delta, x'\in\partial D\cup D_{2r_0},\mathbf{n'}\in N_{\pi(x')}\}
\ee
\be\label{normals2}
N^+_x(q)=\lim_{\delta\downarrow 0}\sup\{\lg q, -\mathbf{n'} \rg: |x'-x|<\delta, x'\in\partial D\cup D_{2r_0},\mathbf{n'}\in N_{\pi(x')}\}.
\ee

We are now ready to give the definition of viscosity solutions.

\bd\label{visdef}
(1) A function $u\in USC((0,T]\times\bard)$ is a viscosity subsolution of Eq.(\ref{pde}) if $u(T,x)\leq h(x)$ for $x\in\bard$, and for any $\phi\in\mathcal{C}^{1,2}((0,T)\times\bard)$, whenever $(t,x)\in(0,T)\times\bard$ is a local maximum of $u-\phi$, then
\ce
-\frac{\partial\phi}{\partial t}(t,x)+H(t,x,u(t,x), D_x\phi, D_x^2\phi)&\leq& 0,~~~\mbox{if}~~~(t,x)\in(0,T)\times D,\\
\min\big\{-\frac{\partial\phi}{\partial t}(t,x)+H(t,x,u(t,x), D_x\phi, D_x^2\phi),N^-_x(D_x\phi)\big\}&\leq& 0, ~~~\mbox{if}~~~(t,x)\in(0,T)\times \partial D.\de

(2) A function $u\in LSC((0,T]\times\bard)$ is a viscosity supersolution of Eq.(\ref{pde}) if $u(T,x)\geq h(x)$ for $x\in\bard$, and for any $\phi\in\mathcal{C}^{1,2}((0,T)\times\bard)$, whenever $(t,x)\in(0,T)\times\bard$ is a local minimum of $u-\phi$, then
\ce
-\frac{\partial\phi}{\partial t}(t,x)+H(t,x,u(t,x), D_x\phi, D_x^2\phi)&\geq&0,~~~\mbox{if}~~~(t,x)\in(0,T)\times D,\\
\max\{-\frac{\partial\phi}{\partial t}(t,x)+H(t,x,u(t,x),D_x\phi, D_x^2\phi),N^+_x(D_x\phi)\big\}&\geq& 0, ~~~\mbox{if}~~~(t,x)\in(0,T)\times \partial D.
\de

(3) A function $u\in \mathcal{C}((0,T]\times\bard)$ is a viscosity solution of  Eq.(\ref{pde}) if and only if it is simultaneously a viscosity sub- and super solution.
\ed

Since below we shall use ``viscosity solution" exclusively,
we shall most time drop the adjunct word ``viscosity".

\br
The above definition coincides essentially with the existing ones for convex domains (see \cite{za12, ren-wuNA, barles2}) and smooth domains
(see \cite{ishiilions, users}).
\er

\br
One can prove, using an argument similar to the proof of \cite[Prop. 2.1]{pa}, that if $H$ is uniformly elliptic, and if $u$ is a subsolution,
then $N_x^-(D_xu)\leq 0$ on the boundary in the viscosity sense. A similar result holds for supersolutions.
\er

\br
It is worth mentioning that oblique derivative problem is investigated by Dupuis and Ishii in \cite{duish1, duish2}. In \cite{duish1} the domain is Lipschitz but the
derivative condition is imposed only in one direction at each point and the direction is assumed to vary smoothly; while in \cite{duish2}
the domain is piecewise smooth and the derivative condition is imposed in several directions at the turning points. Hence neither
of them covers the case of domains considered here.
\er

Then we have the following parabolic analogue to \cite[Prop. 2.6]{koi}. The proof is omitted since it is an obvious
modification of that proposition.

\bp
Let $u\in USC((0,T]\times\bard)$ (resp. $u\in LSC((0,T]\times\bard)$). Then the following claims are equivalent.

(i) $u$ is a subsolution (resp. supersolution).

(ii) $u(T,\cdot)\leq h(\cdot)$ (resp. $u(T,\cdot)\geq h(\cdot)$) and for $(t,x,p,q,X)\in (0,T)\times \bard\times J^{1,2,+}u(t,x)$ (resp. $J^{1,2,-}u(t,x)$),
\ce
-p+H(t,x,u(t,x),q,X)&\leq& 0 ~(\mbox{resp.} \geq 0),\quad\mbox{if}~~~x\in D,\\
\min\big\{-p+H(t,x,u(t,x),q,X),N_x^-(q)\big\}&\leq&0 \\
~~(\mbox{resp.} \max\{-p+H(t,x,u(t,x),q,X),N_x^+(q)\}&\geq& 0)\quad\mbox{if}~~~x\in \partial D.
\de

(iii) $u(T,\cdot)\leq h(\cdot)$ (resp. $u(T,\cdot)\geq h(\cdot)$) and for $(t,x,p,q,X)\in (0,T)\times \bard\times \bar{J}^{1,2,+}u(t,x)$ (resp. $\bar{J}^{1,2,-}u(t,x)$),
\ce
-p+H(t,x,u(t,x),q,X)&\leq& 0 ~(\mbox{resp.} \geq 0),\quad\mbox{if}~~~x\in D,\\
\min\big\{-p+H(t,x,u(t,x),q,X),N_x^-(q)\big\}&\leq&0 \\
~~(\mbox{resp.} \max\{-p+H(t,x,u(t,x),q,X),N_x^+(q)\}&\geq& 0)\quad\mbox{if}~~~x\in \partial D.
\de
\ep

Denote by $X_n^{t,x,\iota}$ the solution to the penalized SDE:
\ce
X_n^{t,x,\iota}(s)=x+\int_t^sb(r,X_n^{t,x,\iota}(r),\iota(r))dr+\int_t^s\sigma(r,X_n^{t,x,\iota}(r),\iota(r))dw(r)-\frac {n}{2}\int_t^s\nabla\varphi(X_n(r))dr.
\de

Denote by $(X^{t,x,\iota},\xi^{t,x,\iota})$ and $(X^{t',x',\iota},\xi^{t',x',\iota})$ the solutions of (\ref{eq1}), starting from $x$ and $x'$ at time $t$ and $t'$ respectively.

\bp\label{continuity1}
For any $p\geq1$, there exists a constant $C=C_p>0$ such that
\be\label{pmoment}
\sup_{\iota\in\cA_\nu} E\sup_{t\leq s\leq T}|X^{t,x,\iota}(s)|^{2p}<C(1+|x|^{2p}),
\ee
\be\label{ppmoment}
\sup_n\sup_{\iota\in\cA_\nu}E\sup_{s\in[t,T]}|X_n^{t,x,\iota}(s)|^{2p}\leq C_p(1+|x|^{2p}),
\ee
and
\be\label{diffp}
\sup_{\iota\in\cA_\nu} E\sup_{t\vee t'\leq s\leq T}|X^{t,x,\iota}(s)-X^{t',x',\iota}(s)|^{2p}\leq C(|x-x'|^{2p}+|t-t'|^p).
\ee
Moreover, for any $c>0$,
\be\label{expi}
\sup_{\iota\in\cA_\nu} E e^{c |\xi^{t,x,\iota}|_T^t}<\infty.
\ee
\ep

\begin{proof}
In the proof we will denote $X^{t,x,\iota }$ and $\xi^{t,x,\iota}$ by $X^{t,x}$ and $\xi^{t,x}$ respectively for
notational simplicity.

(\ref{pmoment}) is proved in \cite[Lemma 2.8]{aisa} and (\ref{ppmoment}) is derived from (\ref{pmoment}) and
Theorem \ref{pconvergence}.
Denote
\ce
M^{t,x}_s:=x+\int_t^sb(X^{t,x}(r),\iota(r))dr+\int_t^s\sigma(X^{t,x}(r),\iota(r))dw(r).
\de

Assume for simplicity $t<t'$ and Note that
\ce
X^{t,x}(s)-X^{t',x'}(s)&=&X^{t,x}(t')-x'+\xi^{t,x}(t')+\int_{t'}^s\big(b(X^{t,x}(r),\iota(r))-b(X^{t',x'}(r),\iota(r))\big)dr\\
&&+\int_{t'}^s\big(\sigma(X^{t,x}(r),\iota(r))-\sigma(X^{t',x'}(r),\iota(r))\big)d w(r)-\big(\xi^{t,x}(s)-\xi^{t',x'}(s)\big).
\de
Applying It\^o's formula to $e^{-\frac1{\gamma}[f(X^{t,x}(s))+f(X^{t',x'}(s))]}|X^{t,x}(s)-X^{t',x'}(s)|^{2}$ and using condition (C),
we get for any $p\geq1$,
\ce
&&E\sup_{t\vee t'\leq s\leq T}|X^{t,x}(s)-X^{t',x'}(s)|^{2p}\\&\leq&CE|X^{t,x}(t')-x'+\xi^{t,x}(t')|^{2p}\\
&\leq &CE|X^{t,x}(t')-x'|^{2p}+E(|\xi^{t,x}(t')|^{2p})\\
&\leq&C|x-x'|^{2p}e^{C_1|t-t'|^p}+E(|\xi^{t,x}|_{t'}^t)^{2p}\\
&\leq&C(|x-x'|^{2p}+|t-t'|^p),
\de
where the last inequality holds since according to \cite[lemma 2.3]{aisa},
\ce
E(|\xi^{t,x}|_{t'}^t)^{2p}&\leq& CE\Big[(1+[M^{t,x}]_{\a;t,t'}^{C_2}(t'-t))^{2p}e^{4C_3\Delta_{t,t'}(M^{t,x})}\{\Delta_{t,t'}(M^{t,x})\}^{2p}\Big]\\
&\leq&C|t-t'|^p,~~~~\forall\alpha\in(0,1).
\de
This proves (\ref{diffp}). For the last one, notice that by Proposition 2.1 in \cite{ren-wuAP} there exists a constant $\mu>0$ such that
\ce
E\Big(e^{\mu (|\xi^{t,x}|_T^t)^2}\Big)<\infty.
\de
Then for any $c>0$,
\ce
E\Big(e^{c |\xi^{t,x}|_T^t}\Big)\leq E\Big(e^{\mu (|\xi^{t,x}|_T^t)^2+\frac{c^2}{4\mu}}\Big)<\infty.
\de
It is easy to see that all the constants appearing in the proof are independent of $\iota$ and we are done.
\end{proof}

Assume $g$ and $h$ are continuous functions satisfying that for all $(t,x,\iota)\in [0,T]\times\bard\times U$, there exist some $C>0$ and $p\geq1$,
\ce
|g(t,x,\iota)|+|h(x)|\leq C(1+|x|^p).
\de

Let
\ce
 V_n(t,x):=\inf_{\nu}\inf_{\iota\in \mathcal{A}_\nu}E\Big(\int_t^T g(s,X_n^{t,x,\iota}(s),\iota(s))ds+h(X_n^{t,x,\iota}(T))\Big),
\de
and \ce
 V_\nu(t,x):=\inf_{\iota\in \mathcal{A}_\nu}E\Big(\int_t^T g(s,X^{t,x,\iota}(s),\iota(s))ds+h(X^{t,x,\iota}(T))\Big).
\de
Then $V=\inf_\nu V_\nu$ and we have the following dynamic programming principle.

\bp\label{dpp}
The function $V$ satisfies the dynamic programming principle. That is, for every $(t,x)\in(0,T]\times\bard$ and any stopping time $\tau$ valued in $[t,T]$, the following hold:\\
(1) For all $\nu$, $\iota\in\mathcal{A}_\nu$,
\ce
V(t,x)\leq E\Big(\int_t^\tau g(s,X^{t,x,\iota}(s),\iota(s))ds+V(\tau, X^{t,x,\iota}(\tau))\Big).
\de
(2) For any $\delta>0$, there exists $\nu$ and $\exists\iota\in\cA_\nu$ such that
\ce
V(t,x)+\delta\geq E\Big(\int_t^\tau g(s,X^{t,x,\iota}(s),\iota(s))ds+V(\tau, X^{t,x,\iota}(\tau))\Big).
\de
(3) For every $\nu$, $V=V_\nu$.
\ep

\begin{proof}
By Theorem \ref{pconvergence}, $E\sup_{s\in[t,T]}|X^{t,x,\iota}_n(s)-X^{t,x,\iota}(s)|^2\rightarrow0$. Actually, carefully  checking the proof in Section 3, one finds that the convergence holds uniformly with respect to $(t,x,\iota)$ on compact sets:
\be\label{pcon1}
\sup_{(t,x,\iota)\in [0,T]\times B_R(0)\times U}E\sup_{s\in[t,T]}|X_n^{t,x,\iota}(s)-X^{t,x,\iota}(s)|^2\rightarrow0,\quad\mbox{as}~~ n\rightarrow\infty,
\ee
where $B_R(0)$ denotes the closed ball with radius $R$ and center $0$.
 Set \ce\omega'(\delta, R):=\sup_{x,y\in B_R(0), |x-y|\leq\delta,(t,\iota)\in(0,T]\times U}(|g(t,x,\iota)-g(t,y,\iota)|+|h(x)-h(y)|).\de
   Then using (\ref{pmoment}), (\ref{ppmoment}) and (\ref{pcon1}), we get
\ce
|V_n(t,x)-V(t,x)|\leq C\omega'(\delta, R)+C(1+|x|^p)(\frac{C_{n,R}}{\delta}+\frac{1+|x|}{R}),
\de
implying that $V_n\rightarrow V$ uniformly on compact subsets.

Then for any stopping time $\tau\in[t,T]$, by Proposition \ref{continuity1},
\ce
&&E|V_n(\tau,X_n^{t,x,\iota}(\tau))-V(t,X^{\tau,x,\iota}(\tau))|\\
&\leq&\sup_{(s,y)\in(0,T]\times(\bard\cap\bar{B}_R(0))}|V_n(s,y)-V(s,y)|\\
&&+E\Big[|V_n(\tau,X_n^{t,x,\iota}(\tau))-V(\tau,X^{t,x,\iota}(\tau))|\mathbbm{l}_{(|X_n^{t,x,\iota}|_{0,t,T}\vee|X^{t,x,\iota}|_{0,t,T}>R)}\Big]\\
&&+E\Big[|V_n(\tau,X_n^{t,x,\iota}(\tau))-V(\tau,X^{t,x,\iota}(\tau))|\mathbbm{l}_{(|X_n^{t,x,\iota}|_{0,t,T}\vee|X^{t,x,\iota}|_{0,t,T}\leq R)}\Big]\\
&\leq&\sup_{(s,y)\in[0,T]\times(\bard\cap\bar{B}_R(0))}|V_n(s,y)-V(s,y)|\\
&&+C(1+|x|^{p+1})/R+C\omega'(\delta,R)\\
&&+C\Big[E(1+|X^{t,x,\iota}(\tau)|^p+|X_n^{t,x,\iota}(\tau)|^p)^2\Big]^{1/2}\\
&&\cdot\Big(\frac{(E|X^{t,x,\iota}(\tau)-X_n^{t,x,\iota}(\tau)|^2)^{1/2}}{\delta}
+\frac{1+(E|X^{t,x,\iota}(\tau)|^2+E|X_n^{t,x,\iota}(\tau)|^2)^{1/2}}{R}\Big)\\
&\leq&C\omega'(\delta,R)+C(1+|x|^{2p})(1/R+c(n,x)/\delta),
\de
where $c(n,x)\to0$ as $n\to\infty$.

(1) Thus \ce
V(t,x)&\leq&V_n(t,x)+|V_n(t,x)-V(t,x)|\\
&\leq& E\Big(\int_t^\tau g(s,X_n^{t,x,\iota}(s),\iota(s))ds+V_n(\tau, X_n^{t,x,\iota}(\tau))\Big)+|V_n(t,x)-V(t,x)|\\
&\leq& E\Big(\int_t^\tau g(s,X^{t,x,\iota}(s),\iota(s))ds+V(\tau, X^{t,x,\iota}(\tau))\Big)+|V_n(t,x)-V(t,x)|\\
&&+C\omega'(\delta,R)+C(1+|x|^{2p})(1/R+c(n,x)/\delta).
\de
Letting $n\to\infty$ and then $\delta\to0$, $R\to\infty$, we get
\ce
V(t,x)\leq E\Big(\int_t^\tau g(s,X^{t,x,\iota}(s),\iota(s))ds+V(\tau, X^{t,x,\iota}(\tau))\Big).
\de

(2) On the other hand, fix $\delta>0$, then for every $n$, there exists $\nu_n$, and $\exists \iota_n\in\cA_{\nu_n}$ such that
\ce
&&V(t,x)+\delta\geq V_n(t,x)+\delta-|V_n(t,x)-V(t,x)|\\
&\geq&E\Big(\int_t^\tau g(s,X_n^{t,x,\iota_n}(s),\iota_n(s))ds+V_n(\tau, X_n^{t,x,\iota_n}(\tau))\Big)-|V_n(t,x)-V(t,x)|+\delta/2\\
&\geq&E\Big(\int_t^\tau g(s,X^{t,x,\iota_n}(s),\iota_n(s))ds+V(\tau, X^{t,x,\iota_n}(\tau))\Big)-|V_n(t,x)-V(t,x)|+\delta/2\\
&&-C\omega'(\delta,R)-C(1+|x|^{2p})(1/R+c(n,x)/\delta)-\sup_{(s,y)\in(0,T]\times(\bard\cap\bar{B}_R(0))}|V_n(s,y)-V(s,y)|.
\de
Thus we can choose $R, ~\delta'$ and $n$ such that there exist $\nu$, $\iota\in\cA_\nu$ satisfying
\ce
V(t,x)+\delta\geq E\Big(\int_t^\tau g(s,X^{t,x,\iota}(s),\iota(s))ds+V(\tau, X^{t,x,\iota}(\tau))\Big)
\de
and the proof is complete.
\end{proof}

Next we will prove an existence theorem. Set
\be\label{functionh}
H(t,x,r,q,Q):=\max_{\iota\in U}\{-\frac12\mathrm{tr}(\sigma\sigma^*)(t,x,\iota)Q-\lg b(t,x,\iota), q \rg-g(t,x,\iota)\}.
\ee
Note that $V_n$ is a viscosity solution to the following HJB equation:
\be\label{pde-penalized}
\left\{
\begin{array}{lllll}
- \frac{\partial v}{\partial t}+H(t,x,v,Dv,D^2v)+\frac {n}{2}\<\nabla\varphi,Dv\>=0\quad \mbox{in}~~~(0,T)\times \mR^d,\\
\\
v(T,\cdot)=h(\cdot).
\end{array}\right.
\ee

We have
\bt\label{existencethm}

Assume (A), (B), (C) and (H). Then the value function $V$ defined in (\ref{cfV}) is solution to Equation (\ref{pde})
with the above given $H$ .
\et

\begin{proof}
We have by the above argument, that $V_n\rightarrow V$ uniformly on compact subsets.

Thus assume for $\phi\in \mathcal{C}^{1,2}((0,T)\times\mR^d)$, $V-\phi$ attains a a strict local maximum at $(t,x)\in(0,T)\times\overline{D}$.
 If $x\in D$, there exists $(t_n, x_n)\in(0,T)\times \bar{D}$ such that
 $(t_n,x_n)\rightarrow(t,x)$ and $V_n-\phi$ attains a local maximum at $(t_n,x_n)$ and therefore,
 \ce
 &&-\frac{\partial \phi}{\partial t}(t_n,x_n)+H(t_n,x_n, V_n(t_n,x_n),D\phi,D^2\phi)+\frac n 2\<\nabla\varphi(x_n),D_x\phi\>\\
 &=&-\frac{\partial \phi}{\partial t}(t_n,x_n)+H(t_n,x_n, V_n(t_n,x_n),D\phi,D^2\phi)\\
 &\leq&0.
 \de
 Now sending $n\rightarrow\infty$ gives
  \ce
 -\frac{\partial \phi}{\partial t}(t,x)+H(t,x, V, D_x\phi,D_x^2\phi)\leq0,
 \de

  If $x\in \partial D$, we can choose $(t_n, x_n)\in(0,T)\times (D_{1/n}\cup\bard)$ (where as defined in Section 1, $D_{1/n}$ is the $1/n$-neighborhood of $\partial D$) such that
 $(t_n,x_n)\rightarrow(t,x)$ and $V_n-\phi$ attains a local maximum at $(t_n,x_n)$ and therefore,
 \ce
-\frac{\partial \phi}{\partial t}(t_n,x_n)+H(t_n,x_n,V_n,D_x\phi(t_n,x_n),D_x^2\phi(t_n,x_n))+\frac n 2|x_n-\pi(x_n)|\<\mathbf{n}(x_n),D_x\phi(t_n,x_n)\>\leq0,
 \de
 where $\mathbf{n}(y):=\frac{y-\pi(y)}{|y-\pi(y)|}$ for any $y\in D_{3r_0}$.  Letting $n\rightarrow\infty$ gives that for $x\in \partial D$,
 \ce
 \min\{ -\frac{\partial \phi}{\partial t}(t,x)+H(t,x, V, D_x\phi(t,x),D^2_x\phi(t,x)),N^-_x(D\phi(t,x))\}\leq0.
 \de

Thus $V$ is a subsolution to HJB equation (\ref{pde}) in the sense of Definition \ref{visdef}.

Similarly, $V$ is a supersolution to (\ref{pde}) and thus a viscosity solution.
\end{proof}

Now we turn to the uniqueness. To this aim we first establish the following comparison theorem.

\bt\label{comparison}
Suppose $u$ is a subsolution bounded from above and respectively $v$ is  a supersolution bounded from below of (\ref{pde}), then
$u\leq v$ on $(0,T]\times\overline{D}$.
\et

\begin{proof}
Suppose that there exists a $(t_0,x_0)\in (0,T)\times \bard$ such that $u(t_0,x_0)-v(t_0,x_0)>0$. Then there exists a $\delta_0>0$
such that
$$
\theta=u(t_0,x_0)-v(t_0,x_0)-\frac{\delta_0}{t_0}>0.
$$
Set
\ce
\psi(t,s,x,y):=\frac{\alpha}{2}(|x-y|^2+|t-s|^2)-\ve(f(x)+f(y))+\delta_0(\frac{1}{2t}+\frac1{2s}).
\de
We fix $\e\in (0,\frac{\theta}{8})$ temporarily. Since
\ce
\lim_{t\wedge s\downarrow 0}\psi(t,s,x,y)=\infty,
\de
there exists $(t_\a,s_\a,x_\a,y_\a)\in (0,T]\times\bard^2$ such that
\ce
&&u(t_\a,x_\a)-v(s_\a,y_\a)-\psi(t_\a,s_\a,x_\a,y_\a)\\
&=&\max_{(t,s,x,y)\in \big((0,T]\times\bar{D}\big)^2}\big( u(t,x)-v(s,y)-\psi(t,s,x,y)\big)\\
&=&:M_\a.
\de
Of course
\ce
M_\a\geq u(t_0,x_0)-v(t_0,x_0)-\frac{\theta}{4}-\frac{\delta_0}{t_0}=\frac{3}{4}\theta,
\de
which yields that
\be\label{uval}
&&u(t_\a,x_\a)-v(s_\a,y_\a)\no\\&\geq&\frac{3}{4}\theta+\frac{\delta_0}{t_\a}+\frac{\a}{2}(|x_\a-y_\a|^2+|t_\a-s_\a|^2)-\ve(f(x_\alpha)+f(y_\alpha))\no\\
&\geq&\frac12\theta+\delta_0(\frac{1}{2t_\a}+\frac{1}{2s_\a})+\frac{\a}{2}(|x_\a-y_\a|^2+|t_\a-s_\a|^2).
\ee
On the other hand, it is easy to see that $\a\mapsto M_\a$ is decreasing and
$$
M_{\frac\a2}\geq M_\a+\frac\a4(|x_\a-y_\a|^2+|t_\a-s_\a|^2).
$$
Consequently
\be\label{alcon}
\lim_{\a\to\infty}\a(|x_\a-y_\a|^2+|t_\a-s_\a|^2)=0.
\ee
Hence, extracting a subsequence if necessary, we have for some $(\hat{t},\hat{x})\in (0,T]$
\ce
\lim_{\a\to\infty}(t_\a,x_\a)=\lim_{\a\to\infty}(s_\a,y_\a)=(\hat{t},\hat{x}).
\de
It is trivial that $\hat{t}\in (0,T)$, and thus $(\hat{t},\hat{x})\in (0,T)\times\bard$ in fact.

By \cite[Lemma 2.7]{koi} there exist $p,q\in \mR$, $X,Y\in \cS^d$ such that
\ce
(p,\a(x_\a-y_\a)+\ve Df(x_\a), X+\ve D^2f(x_\a))\in \bar{J}^{1,2,+}u(t_\a,x_\a),
\de
\ce
(q,\a(x_\a-y_\a)-\ve Df(y_\a), Y-\ve D^2f(y_\a))\in \bar{J}^{1,2,-}v(s_\a,y_\a),
\de
where
\ce
&&p=\a(t_\a-s_\a)-\delta_0\frac1{2t_\a^2},\\
&&q=\a(t_\a-s_\a)+\delta_0\frac1{2s_\a^2},\\
&&p-q=-\frac{\delta_0}{2}(\frac{1}{t_\a^2}+\frac{1}{s_\a^2}),
\de
and
\ce
&&-3\a\left(
\begin{array}{llll}
{I}\quad 0
\\
0\quad {I}
\end{array}\right)\leq\left(
\begin{array}{llll}
{X}\quad\quad 0
\\
0\quad -{Y}
\end{array}\right)\leq 3\a\left(
\begin{array}{llll}
~~{I}\quad -{I}
\\
-{I}\quad {I}
\end{array}\right) ,
 \de
 If $x_\a\in\partial D$, then for $\mathbf{n}\in N_{x_\a}$, using (\ref{alcon}) we have
 \ce
 &&\<-\mathbf{n}, D_x\psi(t_\a,x_\a,y_\a)\>\\
 &=&\<-\mathbf{n},\a(x_\a-y_\a)-\ve Df(x_\a)\>\\
&\geq &-\frac{\a}{8r_0}|x_\a-y_\a|^2+\frac{\e\gamma}{8r_0}\\
&>& 0
 \de
 if $\a$ is large enough. This implies
 \ce
N^-_{x_\a}(D_x\psi(t_\a,x_\a,y_\a))>0.
\de
Thus
 \be\label{con1}
 -p+H(t_\a,x_\a,\a(x_\a-y_\a)+\ve Df(x_\a), X+\ve D^2f(x_\a))\leq 0.
 \ee
 If $x_\a\in D$, then this inequality holds automatically by definition. Similarly for $\alpha$ large enough,
 \ce
 N^+_{y_\a}(y_\a;-D_y\psi(t_\a,x_\a,y_\a))<0,
 \de
 and then
 \be\label{con2}
 -q+H(t_\a,y_\a,\a(x_\a-y_\a) -\ve Df(y_\a), Y-\ve D^2f(y_\a))\geq 0.
 \ee
Combining the above two inequalities gives
\ce
\frac{\delta_0}{2}(\frac{1}{t_\a^2}+\frac{1}{s_\a^2})&\leq& -H(t_\a,x_\a,u(t_\a,x_\a),\a(x_\a-y_\a)+\ve Df(x_\a), X+\ve D^2f(x_\a))\\
&&+H(t_\a,y_\a,v(t_\a,y_\a),\a(x_\a-y_\a) -\ve Df(y_\a), Y-\ve D^2f(y_\a))\\
&\leq& -H(t_\a,x_\a,u(t_\a, x_\a),\a(x_\a-y_\a),X)+ H(t_\a,y_\a,v(t_\a, y_\a),\a(x_\a-y_\a),Y)+2\om_1(C\ve)\\
&\leq& -\gamma\(u(t_\a,x_\a)-v(t_\a,y_\a)\)+\om_2(|x_\a-y_\a|(1+\a|x_\a-y_\a|))+2\om_1(C\ve),
\de
which together with (\ref{uval}) yields that
\ce
\frac34\gamma\theta+\frac{\delta_0}{2}(\frac{1}{t_\a^2}+\frac{1}{s_\a^2})\leq\om_2(|x_\a-y_\a|(1+\a|x_\a-y_\a|))+2\om_1(C\ve).
\de
Letting $\a\to\infty$ and then $\ve\to0$ we get a contradiction and the proof is complete.
\end{proof}

Now we can state the following

\bt
$V$ is the unique solution.
\et

To prove it, we only need to show that the function given by (\ref{functionh}) verifies the hypotheses (H1)-(H3). But this
is a well known fact, see e.g. \cite[Subsection 3.3.2]{koi}.

\section{Elliptic Equations}

We consider in this section in the same circumstance as in the above one the elliptic equation:
\be\label{pdee}
\left\{
\begin{array}{lllll}
-\lambda u+H(x,u,Du,D^2u)=0\quad \mbox{in}~~~D,\\
\\
-\frac{\partial u}{\partial \mathbf{n}}=0\quad \mbox{on}~~~ \partial D,\\
\end{array}\right.
\ee
where $\lam>0$ and $H(x,r,p,Q)$ is a continuous function from $\bard\times\mR\times \mR^d\times \cS^d$ to $\mR$
satisfying the following hypotheses:

(H4) If $u\leq v$, then there exists a constant $\gamma>0$ satisfying
\be\gamma(v-u)\leq H(x,v,p,X)- H(x,u,p,X),
\ee
for all $x\in \bard$, $p\in \mR^d$, $X\in \cS^d$.

(H5) there exists an increasing continuous function $\omega_1: [0,+\infty)\mapsto [0,+\infty)$, $\omega(0)=0$, such that
\be
|H(x, u,q_1,X)-H(x,u,q_2,Y)|\leq \om_1(|q_1-q_2|+\|X-Y\|)
\ee
for all $(x,u)\in \bard\times\mR$, $p,q\in \mR^d$, $X,Y\in \cS^d$.

(H6) there exists an increasing continuous function $\omega_2: [0,+\infty)\mapsto [0,+\infty)$, $\omega(0)=0$, such that
\be
H(x,u,\a(x-y),Y)-H(y,u,\a(x-y),X)\leq \omega_2(|x-y|(1+\a|x-y|))
\ee
for all $x,y\in D$, $u\in\mR$, $\a>1$ and  $X,Y\in \cS^d$ satisfying
\ce
\label{matrix}
-3\alpha\left(
\begin{array}{llll}
I\quad \quad0
\\
0\quad \quad I
\end{array}\right)\leq  \left(
\begin{array}{llll}
{X}\quad\quad 0
\\
0\quad -{Y}
\end{array}\right)\leq
3\alpha\left(
\begin{array}{llll}
I\quad \quad-I
\\
-I\quad \quad I
\end{array}\right).
\de

For a function $\bard\mapsto \mR$, the elliptic superjet
and subjet are defined respectively by

 \ce
 &&J^{2,+}u(t,x):=\Big\{(q,Q)\in\times\mR^d\times\mathcal{S}^d;\\
~~~~~~~~~~~&&\limsup_{\bard\ni y\to x}\frac{u(y)-u(x)-\lg q,y-x\rg-\frac12\la Q(y-x),y-x\rg}{|x-y|^2}\leq0\Big\},\\
 &&J^{1,2,-}u(t,x):=\Big\{(q,Q)\in\mR\times\mR^d\times\mathcal{S}^d;\\
 ~~~~~~~~~~&&\liminf_{\bard\ni y\to x}\frac{u(y)-u(x)-\lg q,y-x\rg-\frac12\la Q(y-x),y-x\rg}{|x-y|^2}\geq0\Big\}.
 \de
and their closures
\ce
&&\bar{J}^{2,\pm}u(x):=\Big\{(q,Q)\in\mR\times\mR^d\times\mathcal{S}^d;
\exists x_k\in  \bard \mbox{ and } (q_k,X_k)\in  J^{2,\pm}u(x) \\
~~~~~~~~~~&&\mbox{ such that } (x_k, u(x_k),q_k, X_k)\to (x, u(x),q, X)\Big\}.
\de

Parallel, viscosity solutions for (\ref{pdee}) are defined as follows.

\bd\label{visdef}
(1) A function $u\in USC(\bard)$ is called a subsolution of Eq.(\ref{pdee}) if for any $\phi\in\mathcal{C}^{2}(\bard)$,
\ce
-\lam u+H(x,u(x),D_x\phi, D_x^2\phi)&\leq& 0,~~~\mbox{if}~~~x\in D,\\
\min\big\{-\lam u+H(x,u(x),D_x\phi, D_x^2\phi),N^-_x(D_x\phi)\big\}&\leq& 0, ~~~\mbox{if}~~~x\in \partial D
\de
provided $\bard$ is a local maximum of $u-\phi$.

(2) A function $u\in LSC(\bard)$ is called a supersolution of Eq.(\ref{pdee}) if for any $\phi\in\mathcal{C}^{2}(\bard)$,
\ce
-\lam u+H(x,u(x),D_x\phi, D_x^2\phi)&\geq&0,~~~\mbox{if}~~~x\in D,\\
\max\{-\lam u+H(x,u(x),D_x\phi, D_x^2\phi),N^+_x(D_x\phi)\big\}&\geq& 0, ~~~\mbox{if}~~~x\in \partial D
\de
provided  whenever $x\in \bard$ is a local minimum of $u-\phi$, then
(3) A function $u\in \mathcal{C}(\bard)$ is a viscosity solution of  Eq.(\ref{pde}) if and only if it is simultaneously a viscosity sub- and super solution.
\ed

Assume $g(x,\iota)$ is a  continuous function of $x$ such that for all $\iota\in U$, there exist some $C>0$ and $p\geq1$,
\ce
|g(x,\iota)|\leq C(1+|x|^p).
\de
Set \ce
 V_\nu(x):=\inf_{\iota\in \mathcal{A}_\nu}E\Big(\int_0^\infty e^{-\lambda s}g(X^{x,\iota}(s),\iota(s))ds).
\de
Then we can prove along the way paved in the last section.

\bt
$V$ is the unique solution to (\ref{pdee}).
\et

\section{Maximum principle for Neumann problem}
In \cite{ren-wuAP}, we have proved the following result
 \bt\label{support}
 Suppose $D$ is a bounded domain satisfying conditions (A)-(C), $\sigma, ~b\in \mathcal{C}_b^2$, and $(X^{x}, \xi^x)$ is the solution to reflected SDE:
 \be\label{eq1plus}
\begin{cases}
dX(t)=\s(X(t))dw(t)+b(X(t))dt+\xi(t),\\
X(0)=x\in \bard.
\end{cases}
\ee
 Then for any fixed $T>0$,
 the support of $P\circ X^{-1}$  in $\cC([0,T],\mR^d)$ coincides with $\overline{\sS}$ where
\ce
\sS:=\{Z^{x,h},~h\in\cC([0,+\infty),\mR^{d_1}), ~h(0)=0, ~t\rightarrow h(t)~~\mbox{is smooth}\}
\de
and  $(Z^{x,h}, \kappa^{x,h})$ is the solution to the following deterministic Skorohod problem in $D$:
\be\label{skorohodZ}
Z^{x,h}(t)=x+\int_0^t \sigma(Z^{x,h}(s))\dot{h}_sd s+\int_0^t\tilde{b}(Z^{x,h}(s))d s+\kappa^{x,h}(t),
\ee
where $\tilde{b}(x)=b(x)-\frac12 \mathrm{tr}(\sigma\sigma^*(x))$.\et

We will apply this result to establish the maximum principle for PDEs with Neumann conditions. To this aim, we need the following proposition:
\bp\label{submart}
Suppose $X^{x}$ is the solution to the reflected SDE in $D$ with initial state $x$ at time $0$, and $u\in\mathcal{C}([0,\infty)\times\bar{D})$ is a subsolution to the following equation
\be\label{maximumpde2}
-\frac{\partial u}{\partial t}+Lu=0, \quad\mbox{in}\quad [0,\infty)\times D, \quad
-\frac{\partial u}{\partial \mathbf{n}}=0\quad\mbox{on}\quad \partial D.
\ee
Then $t\rightarrow u(t,X^{x}(t))$ is a submartingale. Here  $L:=-\frac12\mathrm{tr}(\sigma\sigma^*)(x)\frac{\partial^2}{\partial x^2}-b(x)\frac{\partial}{\partial x}$. In particular, if
$u\in\mathcal{C}(\bar{D})$ is a subsolution to the following PDE with Neumann condition:
\be\label{maximumpde}
Lu=0, \quad\mbox{in}\quad D, \quad -\frac{\partial u}{\partial \mathbf{n}}=0\quad\mbox{on}\quad \partial D.
\ee
then $t\mapsto u(X^x(t)$ is a submartingale.
\ep

\begin{proof}
Suppose $t>0$ and set for $s\in [0,t]$
$$
v(s,x):=E[u(t,X^{s,x}(t))].
$$
Then $v$ is the solution to (\ref{maximumpde2}) on $[0,t]\times \bar{D}$ with $v(t,x)=u(t,x)$. Hence
$$
v(s,\cdot)\geq u(s,\cdot),~~\forall s\in [0,t].
$$
Therefore by Markov property we have for $s\in [0,t]$,
\ce
E[u(t,X^x(t))|\sF_s]&=&E[u(t, X^{s,X^x(s)}(t))|\sF_s]\\
&=&E[u(t, X^{s,y}(t)]_{y=X^x(s)}\\
&=&v(s,X^x(s))\geq u(s, X^x(s)).
\de
This implies that $t\to u(t,X^x(t))$ is a submartingale.
\end{proof}

\bt\label{smp1}
Let $D$ be an open bounded domain in $\mR^d$ satisfying conditions (A)-(C) and $u\in\mathcal{C}(\bar{D})$ be a viscosity solution to the following Neumann problem:
\be\label{smpeq}
\left\{
\begin{array}{lllll}
Lu=0,  ~~~in ~~D\\
-\frac{\partial u}{\partial \mathbf{n}}=0, ~~~~on ~~~~\partial D.
\end{array}\right.
\ee
If there exists $x_0\in \bar{D}$ at which $u$ attains its maximum in $\bar{D}$. Then
\ce
u(x)=u(x_0)\quad \mbox{for all}\quad x\in {D}_0,
\de
where ${D}_0$ is the closure of the set
\ce &&\big\{y; ~~\exists h\in \cC([0,\infty);\mR^{d_1}) ~~\mbox{smooth}, ~~\exists t_1\geq 0 ~~\mbox{such that} ~y=Z^{x_0,h}(t_1)\\
~&&\mbox{and} ~~ \{Z^{x_0,h}(s), ~s\in[0,t_1]\}\subset\bar{D}\big\}\cap \bar{D}.
\de
\et

\begin{proof}
By Proposition \ref{submart}, $t\rightarrow u(X^{x_0}(t))$ is a submartingale and thus for any $t\geq0$ and any stopping time $\tau$,
\be\label{subineq}
E[u(X^{x_0}(t\wedge \tau))]-u(0,x_0)\geq0.
\ee

Suppose there exists $y\in D_0$ such that
\ce
u(y)=u(Z^{x_0,h}(t_1))<u(x_0).
\de
 Let $\Lambda_y$ be any $\varepsilon$-neighborhood of $Z^{x_0,h}(t_1)$ and set
 \ce
 \tau_y:=\inf\{t\geq0, X^{x_0}(t)\in\Lambda_y\cap\bar{D}\}.
 \de
 By Theorem \ref{support},
 \ce
 P^{x_0}(X^{x_0}(\tau_y\wedge t_1)\in D_0)=1.
 \de
 Note that we can choose an $\varepsilon>0$ such that
 \ce
 u( X^{x_0}(t_1\wedge \tau_y))&=&u(X^{x_0}(t_1\wedge \tau_y))-u(Z^{x_0,h}(t_1))+u(Z^{x_0,h}(t_1))\\
 &<&-\varepsilon+u(x_0).
 \de
 By applying the above two estimates we get
 \ce
&&E^{x_0}[u(X^{x_0}(t_1\wedge \tau_y))-u(x_0)]\\
&\leq&-\varepsilon P^{x_0}(|X^{x_0}(t_1\wedge\tau_y)-Z^{x_0,h}(t_1)|<\ve)\\
&<&0,
\de
which is a contraction to (\ref{subineq}).
\end{proof}

Similarly, we can prove
\bt\label{smp2}
Let $D$ be an open bounded domain in $\mR^d$ satisfying conditions (A)-(C) and $u\in\mathcal{C}([0,\infty)\times\bar{D})$ be a viscosity solution to the following Neumann problem:
\be\label{smpeq}
\left\{
\begin{array}{lllll}
-\frac{\partial u}{\partial t}+Lu=0  ~~~in ~~[0,\infty)\times D\\\\
-\frac{\partial u}{\partial \mathbf{n}}=0,   ~~~~on ~~~~\partial D.
\end{array}\right.
\ee
If there exists $(t_0,x_0)\in [0,\infty)\times\bar{D}$ at which $u$ attains its maximum. Then
\ce
u(t,x)=u(t_0,x_0)\quad \mbox{for all}\quad (t,x)\in {D}_0,
\de
where ${D}_0$ is the closure of the set
\ce &&\Big\{y; ~~\exists h\in \cC([0,\infty);\mR^{d_1}) ~~\mbox{smooth}, ~~\exists t_1\in [t_0,\infty) ~~\mbox{such that} ~y=(t_1, Z^{x_0,h}(t_1-t_0))\\
~&&\mbox{and} ~~\{Z^{x_0,h}(s-t_0), ~s\in[t_0,t_1]\}\subset\bar{D}\Big\}\cap [t_0,\infty)\times\bar{D}.
\de
\et

\begin{proof}
Assume for simplicity $t_0=0$. Otherwise we can replace $u(t,x)$ with $u(t-t_0,x)$.

By Proposition \ref{submart}, $t\rightarrow u(t,X^{x_0}(t))$ is a submartingale and thus for any $t\geq0$ and any stopping time $\tau$,
\be\label{subineq2}
E[u(t\wedge \tau,X^{x_0}(t\wedge \tau))]-u(0,x_0)\geq0.
\ee

Suppose $\exists y\in D_0$, such that
\ce
u(y)=u(t_1,Z^{x_0,h}(t_1))<u(0,x_0).
\de
 Let $\Lambda_y$ be any $\varepsilon$-neighborhood of $Z^{x_0,h}(t_1)$ and set
 \ce
 \tau_y:=\inf\{t\geq0, X^{x_0}(t)\in\Lambda_y\cap\bar{D}\}.
 \de
 By Theorem \ref{support},
 \ce
 P^{x_0}\Big((\tau_y\wedge t_1, X^{x_0}(\tau_y\wedge t_1))\in D_0\Big)=1.
 \de
 And note that for $X^{x_0}(t_1)\in \Lambda_y$,
 \ce
 u(t_1, X^{x_0}(t_1))&=&u(t_1, X^{x_0}(t_1))-u(t_1,Z^{x_0,h}(t_1))+u(t_1,Z^{x_0,h}(t_1))\\
 &<&-\varepsilon+u(0,x_0).
 \de
 By applying the above two estimates we get
 \ce
&&E^{x_0}[u(t_1\wedge \tau_y,X^{x_0}(t_1\wedge \tau))-u(0,x_0)]\\
&<&-\varepsilon P^{x_0}(|X^{x_0}(t_1\wedge\tau_y)-Z^{x_0,h}(t_1)|<\ve)\\
&<&0,
\de
which is a contraction to (\ref{subineq2}).
\end{proof}

\end{document}